\documentclass[11pt]{article}

\usepackage{amssymb}
\usepackage{graphicx}
\usepackage{overcite}

\usepackage{blindtext}
\usepackage{algorithm} 
\usepackage{algorithmic} 
\usepackage{amsmath}

\newtheorem{theo}{Theorem}[section]

\newtheorem{lemm}[theo]{Lemma}

\newtheorem{problem}[theo]{Problem}

\begin{document}

\newenvironment{pf}{
\par
\noindent {\bf Proof.}\rm}%
{\mbox{}\hfill\rule{0.5em}{0.8em} \par \bigskip}

\title{On $t$-relaxed 2-distant circular coloring of graphs\footnote{
Research supported by NSFC 11701080 and 11771080.}}

\author{Dan He\footnote{E-mail: hedanmath@163.com} and Wensong
Lin\footnote{E-mail: wslin@seu.edu.cn}\\
\small School of Mathematics, Southeast University, Nanjing 210096,
P.R. China}

\date{}
\maketitle

\vspace*{-1cm} \setlength\baselineskip{6mm}
\bigskip

\begin{abstract}
Let $k$ be a positive integer. For any two integers $i$ and $j$ in
$\{0,1,\dots,k-1\}$, let $|i-j|_k$ be the circular distance between
$i$ and $j$, which is defined as $\min\{|i-j|,k-|i-j|\}$. Suppose
$f$ is a mapping from $V(G)$ to $\{0,1,\dots,k-1\}$. If, for any two
adjacent vertices $u$ and $v$ in $V(G)$, $|f(u)-f(v)|_k\geq 2$, then
$f$ is called a $\frac{k}{2}$-coloring of $G$. Klostermeyer in
[\citeonline{K1,K2}] introduced a relaxation of
$\frac{k}{2}$-coloring named defective circular coloring. Let $d$
and $t$ be two nonnegative integers. If each vertex $v$ is adjacent
to at most $d$ vertices $u$ with $|f(u)-f(v)|_k\leq 1$, then $f$ is
called a $d$-defective 2-distant circular $k$-coloring, or simply a
$(\frac{k}{2},d)$-coloring of $G$. If $G$ has a
$(\frac{k}{2},d)$-coloring, then we say $G$ is
$(\frac{k}{2},d)$-colorable. In this paper, we give a new relaxation
of $\frac{k}{2}$-coloring. If adjacent vertices receive different
integers, and for each vertex $u$ of $G$, the number of neighbors
$v$ of $u$ with $|f(u)-f(v)|_k=1$ is at most $t$, then $f$ is called
a $t$-relaxed 2-distant circular $k$-coloring, or simply a
$(\frac{k}{2},t)^*$-coloring of $G$. If $G$ has a
$(\frac{k}{2},t)^*$-coloring, then $G$ is called
$(\frac{k}{2},t)^*$-colorable. The minimum integer $k$ such that $G$
is $(\frac{k}{2},t)^*$-colorable is called the $t$-relaxed 2-distant
circular chromatic number of $G$, denoted by $c\chi_2^t(G)$.

In this paper, we determine the $t$-relaxed 2-distant circular
chromatic numbers of paths, cycles and complete graphs. We prove
that, for any two fixed integers $k$ and $t$ with $k\geq2$ and
$t\geq1$, deciding whether $c\chi_2^t(G)\leq k$ for a graph $G$ is
NP-complete expect the case $k=2$ and the case $k=3$ and $t\leq3$,
which are polynomially solvable. For any outerplanar graph $G$, it
is easy to see that $G$ is  $(\frac{6}{2},0)^*$-colorable. We show
that all outerplanar graphs are $(\frac{5}{2},4)^*$-colorable. We
also prove that there is no fixed positive integer $t$ such that all
outerplanar graphs are $(\frac{4}{2},t)^*$-colorable. To compare the
$d$-defective 2-distant circular coloring and the $t$-relaxed
2-distant circular coloring, we also prove that, for any two fixed
positive integers $k$ and $d$, deciding whether a graph $G$ is
$(\frac{k}{2},d)$-colorable is NP-complete expect the case $k\leq3$,
which is polynomially solvable. Moreover, for any outerplanar graph
$G$,
we prove that $G$ is $(\frac{4}{2},2)$-colorable. \\


\end{abstract}

\section{Introduction}

In this paper, we focus on undirected and simple graphs, and we use
standard notations in graph theory (cf. [\citeonline{B-M}]). The
circular chromatic number of a graph is a natural generalization of
the chromatic number of a graph. In [\citeonline{Zhu}], Zhu gave the
definition of circular coloring by introducing the problem of
designing a traffic control system. At a road intersection, each
traffic flow is assigned an interval of time during which it faces a
green light. A complete traffic period is a time interval during
which each traffic flow gets a turn of green light. It needs to
design a red-green light pattern for a complete traffic period and
the pattern will be repeated forever. Assume that each interval of
green light has unit length. The problem of designing a traffic
control system is to minimize the total length of a complete traffic
period.

The problem of designing a traffic control system can be modeled as
a kind of graph coloring problem. We use vertices of a graph to
denote the traffic flow, and two vertices are adjacent if the
corresponding traffic flows are not compatible, which means their
green light intervals must not overlap. A complete traffic period
can be viewed as a circle $C$, and each vertex is assigned an open
arc of $C$ with unit length, such that the intervals assigned to
adjacent vertices are disjoint. The objective is to minimize the
total length of the circle $C$. Let $C$ be a circle of length $r$,
an \textit{$r$-circular coloring} of a graph $G$ is a mapping $f$
which assigns to each vertex $u$ of $G$ an open unit length arc of
$C$, such that for every edge $uv$ of $G$, $f(u)\cap
f(v)=\emptyset$. A graph $G$ is $r$-circular colorable if there is
an $r$-circular coloring of $G$, and the \textit{circular chromatic
number} of $G$, denoted by $\chi_c(G)$, is defined as the minimum
$r$ such that $G$ is $r$-circular colorable.

The circular chromatic number $\chi_c(G)$ of a graph $G$ was first
introduced as ``the star chromatic number'' by Vince in
[\citeonline{Vince}]. Some other equivalent definitions such as
$\frac{k}{q}$-coloring can be seen in [\citeonline{Zhu}]. Let $k$ be
a positive integer. For any two integers $i$ and $j$ in
$\{0,1,\dots,k-1\}$, the \textit{circular distance} between $i$ and
$j$, denoted by $|i-j|_k$, is defined as $\min\{|i-j|,k-|i-j|\}$.
For two integers $q$ and $k$ such $1\leq q\leq k$, a
\textit{$\frac{k}{q}$-coloring} of a graph $G$ is a mapping $f$ from
the vertex set of $G$ to colors in $\{0,1,\dots,k-1\}$, such that
for every edge $uv$ of $G$, $|f(u)-f(v)|_k\geq q$. The
\textit{circular chromatic number} $\chi_c(G)$ is the minimum number
$\frac{k}{q}$ such that $G$ has a $\frac{k}{q}$-coloring. It is easy
to see that a $\frac{k}{1}$-coloring of a graph $G$ is just an
ordinary $k$-coloring of $G$. For any finite graph $G$, the
relationship between its chromatic number and circular chromatic
number is $\chi(G)-1<\chi_c(G)\leq \chi(G)$, which can be found in a
survey of circular coloring by Zhu [\citeonline{Zhu}].

Let $G$ be a graph and $f$ a $\frac{k}{q}$-coloring of $G$. Under
the coloring $f$, for any edge $uv$ of $G$, since the circular
distance between $f(u)$ and $f(v)$ is at least $q$ in the circle of
length $k$, we call coloring $f$ a \textit{$q$-distant circular
$k$-coloring} of $G$. If $G$ has a $\frac{k}{q}$-coloring, we say
$G$ is \textit{$\frac{k}{q}$-colorable}. The minimum integer $k$
such that $G$ is $\frac{k}{q}$-colorable is called the
\textit{$q$-distant circular chromatic number} of $G$, denoted by
$c\chi_q(G)$. In this paper, we focus on the case $q=2$.

In the literature, there are many kinds of relaxation about various
graph colorings, please see [\citeonline{C-W-Z, H-L, L1, L2, L-Z,
W-W}]. Let $k$ and $d$ be two nonnegative integers. In
[\citeonline{C-C-W},\citeonline{C-G-J}], the authors introduced a
relaxed coloring which is named as $(k,d)$-defective coloring. A
\textit{$(k,d)$-defective coloring} of a graph is an assignment of
$k$ colors to the vertices such that each vertex $v$ is adjacent to
at most $d$ vertices having the same color as $v$. Defective
coloring is sometimes known as ``improper coloring'', which is
defined in the book by Jensen and Toft [\citeonline{J-T}]. The
$(k,d)$-defective coloring has been studied extensively in the
literature, please see
[\citeonline{A},\citeonline{F-H},\citeonline{H-K-S}].

There is also a kind of relaxation about circular coloring, which
was called defective circular coloring by Klostermeyer in
[\citeonline{K1,K2}]. \textit{A defective circular coloring} for a
simple graph $G$ is a mapping $f: V(G)\rightarrow\{0,1,\dots,k-1\}$,
such that each vertex $v$ is adjacent to at most $d$ vertices $u$
where $|f(u)-f(v)|_k\geq q$ does not hold. In this paper, we call
this mapping $f$ a \textit{$d$-defective $q$-distant circular
$k$-coloring}, or simply a \textit{$(\frac{k}{q},d)$-coloring}. If
$G$ has a $(\frac{k}{q},d)$-coloring, we say $G$ is
\textit{$(\frac{k}{q},d)$-colorable}. In [\citeonline{K2}],
Klostermeyer gave some fundamental properties of defective circular
coloring and studied the $(\frac{5}{2},d)$-colorability of planar
graphs, $(\frac{5}{2},1)$-colorability of outerplanar graphs and the
$(\frac{5}{2},2)$-colorability of series-parallel graphs. The other
results about defective circular coloring can be seen in
[\citeonline{G-X,S}]. In this paper, we introduce a new relaxation
of circular coloring.

Let $k,q$ and $t$ be nonnegative integers. Suppose $f$ is a mapping
from $V(G)$ to $\{0,1,\dots,k-1\}$. If adjacent vertices receive
different integers, and for each vertex $u$ of $G$, the number of
neighbors $v$ of $u$ with $|f(u)-f(v)|_k<q$ is at most $t$, we say
$f$ is a \textit{$t$-relaxed $q$-distant circular $k$-coloring} of
$G$, or simply a \textit{$(\frac{k}{q},t)^*$-coloring} of $G$. If
$t=0$, the above coloring $f$ is just a $\frac{k}{q}$-coloring of
$G$. In other words, a $(\frac{k}{q},0)^*$-coloring is also a
$\frac{k}{q}$-coloring. If $G$ has a $(\frac{k}{q},t)^*$-coloring,
we say $G$ is \textit{$(\frac{k}{q},t)^*$-colorable}. The minimum
integer $k$ such that $G$ is $(\frac{k}{q},t)^*$-colorable is called
the \textit{$t$-relaxed $q$-distant circular chromatic number} of
$G$, denoted by $c\chi_q^t(G)$. In this paper, we focus on the case
$q=2$. According to the above definition of $t$-relaxed $q$-distant
circular chromatic number, for any finite graph $G$ and positive
integer $t$, we have $\chi(G)\leq c\chi_2^t(G)\leq 2\chi(G)$.
Moreover, $c\chi_2^t(G)=\chi(G)$ when $t\geq\Delta(G)$. If $G$ is
$(\frac{k}{q},t)^*$-colorable, then it is also
$(\frac{k}{q},t)$-colorable, but not vice versa.

Note that the $t$-relaxed $q$-distant circular coloring is different
from the $d$-defective $q$-distant circular coloring. In a
$d$-defective $q$-distant circular coloring, the colors assigned to
two adjacent vertices can be allowed to be the same. However, in
some practical problems, for two adjacent vertices $u$ and $v$,
although the circular distance between colors assigned to $u$ and
$v$ may not be required to be greater than or equal to $q$, it is
often necessary to require that $u$ and $v$ must receive different
colors. In a $t$-relaxed $q$-distant circular coloring, we put
forward this requirement. This means the constraint in $t$-relaxed
$q$-distant circular coloring is stronger than that in $d$-defective
2-distant circular coloring. In section 5, we will give an example
of a graph which is $(\frac{4}{2},1)$-colorable, but is not
$(\frac{4}{2},1)^*$-colorable.

Suppose $f$ is a $t$-relaxed (resp. $d$-defective) 2-distant
circular $k$-coloring of a graph $G$. Let $v$ be a vertex of $G$. If
$v$ has $p$ neighbors $u$ with $|f(u)-f(v)|_k=1$ (resp.
$|f(u)-f(v)|_k\leq1$), we say $v$ is relaxed $p$ times, or $v$ has
$p$ relaxations. Assume $S$ is a subset of $V(G)$. If there are $p$
vertices $w\in N(v)\cap S$ with $|f(w)-f(v)|_k=1$ (resp.
$|f(w)-f(v)|_k\leq1$), we say $v$ has $p$ relaxations in $S$. If,
for each neighbor $u$ of $v$, $|f(u)-f(v)|_k\geq2$, then the vertex
$v$ is not relaxed. For a mapping from the vertex set to
$\{0,1,\dots,k-1\}$, if there exists at least one vertex $u$ which
has relaxations, this mapping is called a \textit{relaxed 2-distant
circular $k$-coloring}.

We proceed as follows. In section 2 we determine $t$-relaxed
2-distant circular chromatic numbers of paths, cycles and complete
graphs. In section 3 we study the complexity of $t$-relaxed
2-distant circular coloring problem. We prove that for two fixed
integers $k\geq2$ and $t\geq1$, the problem of deciding whether a
graph $G$ is $(\frac{k}{2},t)^*$-colorable is NP-complete expect the
case $k=2$ and the case $k=3$ and $t\leq3$, which are polynomially
solvable. To compare the above two relaxations of circular coloring,
in section 4, we prove that, for any two fixed positive integers $k$
and $d$, deciding whether a graph $G$ is $(\frac{k}{2},d)$-colorable
is NP-complete expect the case $k\leq3$, which is polynomially
solvable. In section 5 we investigate relaxed 2-distant circular
colorings of outerplanar graphs. For any outerplanar graph $G$, it
is easy to see $c\chi_2(G)\leq6$. This means outerplannar graphs are
all $(\frac{6}{2},0)^*$-colorable. When $k=5$, we give an
outerplanar graph which is not $(\frac{5}{2},1)^*$-colorable. With
the aid of ordered breath first search technique, we construct a
polynomial-time algorithm to obtain a $(\frac{5}{2},4)^*$-coloring
for any outerplanar graph. When $k=4$, for any fixed positive
integer $t$, we present a class of outerplanar graphs that is not
$(\frac{4}{2},t)^*$-colorable. However, on the other hand, we prove
that all outerplanar graphs are $(\frac{4}{2},2)$-colorable.

\section{$c\chi_2^t(G)$ for graph $G$ in some special classes of graphs}

In this section, we determine the $t$-relaxed 2-distant circular
chromatic numbers of paths, cycles and complete graphs.

\begin{theo}\label{Pn}
Let $t$ be a positive integer and $P_n$ a path on $n$ vertices with
$n\geq3$. Then
$$c\chi_2^t(P_n)=\left\{
                            \begin{array}{lll}
                              4, &  \mbox{if}~~ t=1,\\
                              2, &  \mbox{if}~~ t\geq2.
                           \end{array}
                          \right.$$
\end{theo}

\begin{pf} Since the chromatic number of path $P_n$ is 2, we have $2\leq c\chi_2^t(P_n)\leq
4$. First we deal with the case $t=1$. Suppose to the contrary that
$c\chi_2^1(P_n)\leq3$. Let $f: V(P_n)\rightarrow \{0,1,2\}$ be a
$(\frac{3}{2},1)^*$-coloring of $P_n$. Since the circular distance
between any two colors in $\{0,1,2\}$ is equal to 1, the vertex in
$V(P_n)$ with degree 2 can not be colored by colors in $\{0,1,2\}$.
This is a contradiction and so $c\chi_2^1(P_n)=4$. Next we assume
$t\geq2$. By assigning the color 0 to all vertices in one part of
$P_n$ and the color 1 to all vertices in the other part of $P_n$, we
obtain a $(\frac{2}{2},t)^*$-coloring of $P_n$. Thus,
$c\chi_2^t(P_n)=2$ for $t\geq2$.
\end{pf}

\begin{theo}\label{Cn}
Let $t$ be a positive integer and $C_n$ a cycle on $n$ vertices.
Then
$$c\chi_2^1(C_n)=\left\{
                            \begin{array}{lll}
                              4, & \mbox{if}~~ n>3,\\
                              5, &  \mbox{if}~~n=3
                           \end{array}
                          \right. ~~\mbox{and }~~c\chi_2^t(C_n)=\left\{
                            \begin{array}{lll}
                              2, &  \mbox{if}~~t\geq2~ \mbox{and}~ n~ \mbox{is even},\\
                              3, & \mbox{if}~~  t\geq2~ \mbox{and}~ n~ \mbox{is odd}.
                           \end{array}
                          \right.$$
\end{theo}

\begin{pf} Note that the chromatic number of $C_n$ is 2 if $n$ is
even and is 3 if $n$ is odd. Since the degree of every vertex in
$C_n$ is 2, we have $c\chi_2^t(C_n)=\chi(C_n)$ for any $t\geq2$.
Therefore, it is only necessary to consider the case $t=1$. From
Theorem \ref{Pn}, we have $c\chi_2^1(P_n)=4$ for $n\geq3$,
indicating that $c\chi_2^1(C_n)\geq4$ for $n\geq3$. Let
$C_n=u_1u_2\dots u_n$. If $n$ is even (resp. odd and $n\geq5$), then
we color the vertices $u_1,u_2,\dots,u_n$ in the manner
$0,2,0,2\dots,0,2$ (resp. $0,2,0,2,\dots,0,2,0,1,3)$. It is
straightforward to check that this gives a
$(\frac{4}{2},1)^*$-coloring of $C_n$. Thus, $c\chi_2^1(C_n)=4$ if
$n\geq4$. When $n=3$, suppose to the contrary that
$c\chi_2^1(C_3)=4$. Let $f: V(C_3)\rightarrow \{0,1,2,3\}$ be a
$(\frac{4}{2},1)^*$-coloring of $C_n$. Without loss of generality,
assume $f(u_1)=0$. Since each vertex has at most one relaxation
under the coloring $f$, the color 2 must be assigned to $u_2$ or
$u_3$. Then whether $u_2$ or $u_3$ is colored by 2, another vertex
can not be colored by colors in $\{0,1,2,3\}$. This is a
contradiction and $c\chi_2^1(C_3)\geq5$. By assigning colors $0,2$
and $4$ to vertices $u_1,u_2$ and $u_3$ respectively, we get a
$(\frac{5}{2},1)^*$-coloring of $C_3$. Thus, $c\chi_2^1(C_3)=5$ and
the theorem holds.
\end{pf}

\begin{theo}\label{Kn}
Let $t$ be a positive integer and $K_n$ a complete graph on $n$
vertices. Then
$$c\chi_2^t(K_n)=\left\{
                            \begin{array}{lll}
                              \frac{3n}{2}, & \mbox{if}~~t=1~\mbox{and $n$ is
                              even},\\[2mm]
                              \frac{3n+1}{2}, & \mbox{if}~~ t=1~\mbox{and $n$ is
                              odd},\\[2mm]
                              n, & \mbox{if}~~ t\geq2.
                           \end{array}
                          \right.$$
\end{theo}

\begin{pf} Note that each color $c$ in $\{0,1,\dots,n-1\}$ has only two
colors whose circular distance from $c$ is equal to 1. By assigning
different colors in $\{0,1,\dots,n-1\}$ to $n$ vertices of $K_n$, we
get a 2-relaxed 2-distant circular $n$-coloring. It follows that
$c\chi_2^t(K_n)=n$ for $t\geq2$.

Next we deal with the case $t=1$. Suppose $c\chi_2^1(K_n)=k$ and let
$f: V(K_n)\rightarrow \{0,1,\dots,k-1\}$ be a
$(\frac{k}{2},1)^*$-coloring of $K_n$. By the definition of
1-relaxed 2-distant circular coloring, the colors assigned to the
$n$ vertices must be different from each other, and for each vertex
$u\in V(K_n)$, there is at most one neighbor $v$ of $u$ satisfing
$|f(u)-f(v)|_k=1$. Let $j$ be a color in $\{0,1,\dots,k-1\}$. Denote
by $n_j$ the number of vertices assigned color $j$ by $f$. Then it
is easy to see that $n_j+n_{j+1}+n_{j+2}\leq 2$, where the
subscripts $j+1$ and $j+2$ are taken module $k$. Thus
$$3n=\displaystyle\sum\limits_{j=0}^{k-1}(n_j+n_{j+1}+n_{j+2})\leq 2\cdot k,$$
implying $k\geq\frac{3}{2}n$. This means $c\chi_2^1(K_n)\geq
\frac{3}{2}n$. To prove the theorem, it is necessary to give a
1-relaxed 2-distant circular coloring of $K_n$.

If $n$ is even, let $n=2p$ ($p$ is an integer greater than 1) and
$C_1=\{0,1,3,4,\dots,3p-3,3p-2\}$. Note that $|C_1|=2p=n$. By
assigning different colors in $C_1$ to vertices of $K_n$, we get a
$(\frac{3p}{2},1)^*$-coloring of $K_n$ and so
$c\chi_2^1(K_n)=3p=\frac{3n}{2}$ for even $n$. If $n$ is odd, it is
clear that $c\chi_2^1(K_n)=k\geq \frac{3n+1}{2}$. Let $n=2p+1$ ($p$
is a positive integer) and $C_2=\{0,1,3,4,\dots,3p-3,3p-2,3p\}$.
Then $|C_2|=2p+1=n$. By assigning different colors in $C_2$ to
vertices of $K_n$, we get a $(\frac{3p+2}{2},1)^*$-coloring and so
$c\chi_2^1(K_n)=3p+2=\frac{3n+1}{2}$ for odd $n$.
\end{pf}

\section{The complexity of $t$-relaxed 2-distant circular coloring }

The complexity of coloring and relaxed coloring has been studied
extensively in the literature [\citeonline{B-H,C-G-J,G-J,G-J-S}].
Let's now introduce some known results about coloring and relaxed
coloring of graphs, which are quite essential for proving our
results on $t$-relaxed 2-distant circular coloring. Let $l$ be a
fixed integer with $l\geq 3$. It is well known that the $l$-coloring
problem is NP-complete, please see [\citeonline{G-J}]. Furthermore,
Garey etc. proved the following result in [\citeonline{G-J-S}].

\begin{lemm}\label{kc-complexity}  [\citeonline{G-J-S}]
The problem of determining whether a graph is $3$-colorable is
NP-complete, even for planar graphs whose maximum degrees are at
most four.
\end{lemm}

In section 1, we have mentioned the $(k,t)$-defective coloring and
2-distant circular coloring. The complexity about these two kinds of
coloring has been proved as follows.

\begin{lemm} \label{kt-complexity} [\citeonline{C-G-J}]
For any two fixed integers $k\geq2$ and $t\geq1$, the problem of
determining whether a graph is $(k,t)$-defective colorable is
NP-complete.
\end{lemm}

\begin{lemm}\label{comp-k-cir} [\citeonline{B-H}]
Let $G$ be a graph. For any fixed positive integer $k$ with
$k\geq5$, it is NP-complete to decide whether $c\chi_2(G)\leq k$.
\end{lemm}

In this section, we consider the complexity of $t$-relaxed 2-distant
circular coloring. Let $k\geq 2$ and $t\geq1$ be two fixed integers.

\begin{problem}\label{dis-rel-cir-pro}
$((\frac{k}{2},t)^*\mbox{-coloring problem})$

\noindent Instance: A graph $G$.

\noindent Question: Is $G$ $(\frac{k}{2},t)^*$-colorable?
\end{problem}

It is obvious that if $G$ is a nonempty graph, then
$c\chi_2^t(G)\geq2$ for any positive integer $t$. Based on this
observation, we assume $k\geq2$ in this section. When $k=2$ (resp.
$k=3$), the circular distance between any two colors in $\{0,1\}$
(resp. $\{0,1,2\}$) is $1$. Note that in a $t$-relaxed 2-distant
circular coloring, adjacent vertices must receive different colors.
The following lemma is easy to see.

\begin{lemm} \label{rel-k-23}
Let $t$ be a nonnegative integer and $\Delta(G)$ be the maximum
degree of $G$. Then a graph $G$ is $(\frac{2}{2},t)^*$-colorable if
and only if $G$ is a bipartite graph with $\Delta(G)\leq t$. And a
graph $G$ is $(\frac{3}{2},t)^*$-colorable if and only if $G$ is a
3-partite graph with $\Delta(G)\leq t$.
\end{lemm}

Based on Lemma \ref{rel-k-23}, Problem \ref{dis-rel-cir-pro} for
$k=2$ is polynomially solvable, and a graph $G$ is
$(\frac{3}{2},t)^*$-colorable if and only if $G$ is 3-colorable and
$\Delta(G)\leq t$. Note that if $\Delta(G)=3$, then $G$ is
3-colorable unless $G$ is isomorphic to complete graph $K_4$. Thus
when $k=3$ and $t\leq3$, Problem \ref{dis-rel-cir-pro} is
polynomially solvable. From Lemma \ref{kc-complexity}, we know that
3-coloring problem is NP-complete for graphs whose maximum degrees
are at most 4. Combining the above discussion and the fact that
$c\chi_2^t(G)=\chi(G)$ for any $t\geq\Delta(G)$, we have the
following result.

\begin{lemm} \label{rel-k-3}
Let $k$ and $t$ be two integers with $k\geq2$ and $t\geq1$. If $k=2$
or $k=3$ and $t\leq3$, the $(\frac{k}{2},t)^*$-coloring problem is
polynomially solvable. Meanwhile, if $t\geq4$, the
$(\frac{3}{2},t)^*$-coloring problem is NP-complete.
\end{lemm}

From now on, we assume $k\geq4$. We consider the following three
cases: $k=4$ and $t=1$, $k=4$ and $t\geq2$, and $k\geq5$. The
following lemma is essential in our proofs.

\begin{lemm} \label{P4}
Let $P_4=uxyv$ be a path on four vertices. Then the following two
properties hold:

\noindent (1) If $f$ is a $(\frac{4}{2},1)^*$-coloring of $P_4$ with
$f(u)=f(v)$, then each of $u$ and $v$ must have one relaxation.

\noindent (2) For any two different colors $a$ and $b$ in
$\{0,1,2,3\}$, there is a $(\frac{4}{2},1)^*$-coloring of $P_4$ with
$\{f(u),f(v)\}=\{a,b\}$ such that both $u$ and $v$ have no
relaxations.
\end{lemm}

\begin{pf} (1) Suppose $f(u)=f(v)=a$. If $u$ has no relaxations,
then $f(x)=a+2 \pmod4$, implying $f(y)\in \{a+1 \pmod 4,
a+3\pmod4\}$. So $y$ has two relaxations, which contradicts the fact
that $f$ is a $(\frac{4}{2},1)^*$-coloring. It follows that $f(x)\in
\{a+1 \pmod 4, a+3\pmod4\}$. Similarly, $f(y)\in \{a+1 \pmod 4,
a+3\pmod4\}$. Since $f(x)\neq f(y)$, we must have
$\{f(x),f(y)\}=\{a+1 \pmod 4, a+3\pmod4\}$. This means each of $u$
and $v$ has exactly one relaxation.

(2) We only need to consider two cases: $f(u)=a$ and $f(v)=a+1\pmod
4$, or $f(u)=a$ and $f(v)=a+2\pmod 4$. If $f(u)=a$ and
$f(v)=a+1\pmod 4$, let $f(x)=a+2\pmod 4$ and $f(y)=a+3\pmod 4$. If
$f(u)=a$ and $f(v)=a+2\pmod 4$, let $f(x)=a+2\pmod 4$ and $f(y)=a$.
It is easy to see that $f$ is the desired
$(\frac{4}{2},1)^*$-coloring of $P_4$.
\end{pf}

\begin{lemm} \label{k-four-t-1}
The $(\frac{4}{2},1)^*$-coloring problem is NP-complete.
\end{lemm}

\begin{pf} We prove this lemma by reducing the $(4,1)$-defective coloring
problem to the $(\frac{4}{2},1)^*$-coloring problem. Given an
instance $G$ of the $(4,1)$-defective coloring problem, we construct
an instance $G^*$ of the $(\frac{4}{2},1)^*$-coloring problem by
replacing each edge of $G$ with a path $P_4$. It is obvious that the
reduction can be accomplished in polynomial time.

If $G$ has a $(4,1)$-defective coloring $f$ using colors in
$\{0,1,2,3\}$, then we can define a $(\frac{4}{2},1)^*$-coloring $g$
of $G^*$ as follows. For each vertex $v\in V(G)$, let $g(v)=f(v)$.
Suppose the edge $uv$ of $G$ is replaced with the path $P_4=uxyv$.
The vertices $x$ and $y$ can be colored as follows.

If $f(u)=f(v)=a$, then let $g(x)=a+1\pmod4$ and $g(y)=a+3\pmod4$. If
$f(u)\neq f(v)$, then by Lemma \ref{P4} (2), we can color the two
vertices $x$ and $y$ properly such that both $u$ and $v$ have no
relaxations. It is clear that the mapping $g$ is a
$(\frac{4}{2},1)^*$-coloring of $G^*$.

On the other hand, suppose $G^*$ has a $(\frac{4}{2},1)^*$-coloring
$g$, then let $h(u)=g(u)$ for any vertex $u\in V(G)$. Let $u$ and
$v$ be any two adjacent vertices in $G$. By Lemma \ref{P4} (1), we
know that if $g(u)=g(v)$, then there is no any other neighbors of
$u$ (resp. $v$) in $G$ which has the same color as $u$ (resp. $v$).
Thus the mapping $h$ is a $(4,1)$-defective coloring of $G$.

We have shown that $G$ has a $(4,1)$-defective coloring if and only
if $G^*$ has a $(\frac{4}{2},1)^*$-coloring. Consequently, the
NP-completeness of the $(\frac{4}{2},1)^*$-coloring problem follows
from that of the $(4,1)$-defective coloring problem. Therefore, the
lemma holds by Lemma \ref{kt-complexity}.
\end{pf}

We now turn to the case $k=4$ and $t\geq2$.

\begin{lemm} \label{k-four-t}
Let $t$ be a fixed integer with $t\geq2$. Then the
$(\frac{4}{2},t)^*$-coloring problem is NP-complete.
\end{lemm}

\begin{pf}
We prove this lemma by the reduction from the $(4,t)$-defective
coloring problem to the $(\frac{4}{2},t)^*$-coloring problem. Let
$G$ be an instance of the $(4,t)$-defective coloring problem. We
construct an instance $G_t^*$ from $G$ by replacing each edge $uv\in
E(G)$ with the graph $A_{t-1}(u,v)$, which consists of a path
$P_4=uxyv$ with $x$ (resp. $y$) connecting to the $2(t-1)$ vertices
of the $t-1$ disjoint copies of $K_2$. Please see Figure \ref{A4uv}
for the graph $A_4(u,v)$ as an example. It is obvious that the
reduction is polynomial.

 \begin{figure}[h!]
 \centering \resizebox{7.2cm}{3.1cm}{\includegraphics{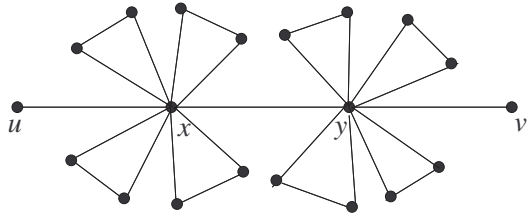}}
 \caption{The graph $A_4(u,v)$}\label{A4uv}
 \end{figure}

Suppose $G$ has a $(4,t)$-defective coloring $f$ using colors in
$\{0,1,2,3\}$. Then we define a $(\frac{4}{2},t)^*$-coloring $g$ of
$G_t^*$ as follows. For each vertex $v$ in $V(G)$, let $g(v)=f(v)$.
For any edge $uv\in E(G)$, there is a path $P_4=uxyv$ in $G_t^*$. If
$f(u)=f(v)=a$, then let $g(x)=a+1\pmod4$ and $g(y)=a+3\pmod4$, and
assign $a+2\pmod4$ and $a+3\pmod4$ (resp. $a$ and $a+1\pmod4$) to
the two vertices in each of the $t-1$ copies of $K_2$ connected to
$x$ (resp. $y$). If $f(u)\neq f(v)$, by Lemma \ref{P4} (2), we can
assign colors to $x$ and $y$ satisfying
$|g(u)-g(x)|_4=|g(v)-g(y)|_4=2$ and $|g(x)-g(y)|_4\geq1$. We then
assign $g(x)+1\pmod4$ and $g(x)+2\pmod4$ (resp.  $g(y)+1\pmod4$ and
$g(y)+2\pmod4$) to the two vertices in each of the $t-1$ copies of
$K_2$ connected to $x$ (resp. $y$). Since $f$ is a $(4,t)$-defective
coloring of graph $G$, it is easy to see that $g$ is a
$(\frac{4}{2},t)^*$-coloring of graph $G_t^*$.

Now suppose $G_t^*$ has a $(\frac{4}{2},t)^*$-coloring $g$. Let $uv$
be an edge of $G$ and $uxyv$ be the corresponding $P_4$ in $G_t^*$.
By Lemma \ref{P4} (1), if $g(u)=g(v)$, then
$|g(u)-g(x)|_4=|g(v)-g(y)|_4=1$. It follows that, for each vertex
$u$ of $G$, there are at most $t$ neighbors $v$ in $G$ with
$g(u)=g(v)$. Let $h(v)=g(v)$ for each $v\in V(G)$. It is clear that
$h$ is a $(4,t)$-defective coloring of $G$.

Since we have shown that $G$ has a $(4,t)$-defective coloring if and
only if $G_t^*$ has a $(\frac{4}{2},t)^*$-coloring, the
NP-completeness of the $(\frac{4}{2},t)^*$-coloring is established
by Lemma \ref{kt-complexity}.
\end{pf}

Finally, we consider the general case $k\geq5$. We construct a
reduction from the 2-distant circular coloring problem to the
$(\frac{k}{2},t)^*$-coloring problem. Let $G$ and $H$ be two simple
graphs. The \textit{composition} of $G$ and $H$ is the graph $G[H]$
with vertex set $V(G)\times V(H)$ in which $(u,v)$ is adjacent to
$(u',v')$ if and only if either $uu'\in E(G)$ or $u=u'$ and $vv'\in
E(H)$. Denote the empty graph on $p$ vertices by $K_p^c$.

\begin{lemm} \label{k-five-t}
Let $k$ and $t$ be two fixed integers with $k\geq5$ and $t\geq1$.
Then the $(\frac{k}{2},t)^*$-coloring problem is NP-complete.
\end{lemm}

\begin{pf}
We prove this lemma by reducing the $\frac{k}{2}$-coloring problem
to the $(\frac{k}{2},t)^*$-coloring problem. Let $G$ be an instance
of the $\frac{k}{2}$-coloring problem. We construct an instance
$G^*=G[K_{kt+1}^c]$ of the $(\frac{k}{2},t)^*$-coloring problem.

Assume $|V(G)|=n$ and let $V(G)=\{v_1,v_2,\dots,v_n\}$. Let the
vertex set of $K_{kt+1}^c$ be $\{w_1,w_2,\dots,w_{kt+1}\}$. For
simplicity, for $i\in\{1,2,\dots,n\}$ and $j\in \{1,2,\dots,kt+1\}$,
let $v_i^j$ denote the vertex $(v_i,w_j)$ of $G[K_{kt+1}^c]$. We use
$V_i$ to denote the vertex set $\{v_i^1,v_i^2,\dots,v_i^{kt+1}\}$.

Suppose $f$ is a $\frac{k}{2}$-coloring coloring of $G$. Then we
define $g(v_i^j)=f(v_i)$ for $i=1,2,\dots,n$ and $j=1,2,\dots,kt+1$.
Since $V_i$ is an independent set of $G^*$, it is obvious that $g$
is a 2-distant circular $k$-coloring of $G^*$. Of course, $g$ is a
$(\frac{k}{2},t)^*$-coloring.

Now suppose $g$ is a $(\frac{k}{2},t)^*$-coloring of $G^*$. We are
going to define a 2-distant circular $k$-coloring $h$ of $G$. For
$i=1,2,\dots,n$, since $V_i$ has $kt+1$ vertices and $g$ only use
$k$ colors, by the pigeonhole principle, there is a color that is
assigned to at least $t+1$ vertices in $V_i$, denoted this color by
$c_i$. Let $h(v_i)=c_i$ for $i=1,2,\dots,n$. If $v_iv_s\in E(G)$, we
have $|h(v_i)-h(v_s)|_k\geq2$, since otherwise the vertex in $V_i$
with color $c_i$ would be relaxed at least $t+1$ times in $V_s$
under the $(\frac{k}{2},t)^*$-coloring $g$ of $G^*$. Thus $h$ is a
2-distant circular $k$-coloring of $G$.

We have shown that $G$ has a $\frac{k}{2}$-coloring if and only
$G^*$ has a $(\frac{k}{2},t)^*$-coloring. So, the NP-completeness of
the $(\frac{k}{2},t)^*$-coloring is establish from Lemma
\ref{comp-k-cir}.
\end{pf}

Combining Lemma \ref{rel-k-3}, \ref{k-four-t-1}, \ref{k-four-t} and
\ref{k-five-t}, we have the following theorem.

\begin{theo} \label{t-relexed-comp}
Let $k$ and $t$ be two fixed integers with $k\geq2$ and $t\geq1$.
The $(\frac{k}{2},t)^*$-coloring problem is NP-complete expect the
case $k=2$ and the case $k=3$ and $t\leq3$, which are polynomially
solvable.
\end{theo}

\section{The complexity of defective circular coloring }

To compare the difference between $t$-relaxed 2-distant circular
coloring and $d$-defective 2-distant circular coloring, we consider
the complexity of $d$-defective 2-distant circular coloring in this
section. Let $k$ and $d$ be two fixed positive integers.

\begin{problem}\label{dis-def-cir-pro}
$((\frac{k}{2},d)\mbox{-coloring problem})$

\noindent Instance: A graph $G$.

\noindent Question: Is $G$ $(\frac{k}{2},d)$-colorable?
\end{problem}

Note that, for $1\leq k\leq 3$, the circular distance between any
two colors in $\{0,1,\dots,k-1\}$ are less than $2$. Therefore, if a
graph $G$ is $(\frac{k}{2},d)$-colorable for some $1\leq k\leq3$,
then we must have $\Delta(G)\leq d$. On the other hand, it is easy
to see that if $\Delta(G)\leq d$, then $G$ is
$(\frac{k}{2},d)$-colorable for any $1\leq k\leq3$. Thus, a graph
$G$ is $(\frac{k}{2},d)$-colorable for some $1\leq k\leq3$ if and
only if $\Delta(G)\leq d$. So, we have the following lemma.

\begin{lemm} \label{def-k-3} Let $k$ and $d$ be two positive
integers. If $k\leq3$, then the $(\frac{k}{2},d)$-coloring problem
is polynomially solvable.
\end{lemm}

\begin{lemm} \label{k-four-d}
Let $d$ be a fixed positive integer. Then the
$(\frac{4}{2},d)$-coloring problem is NP-complete.
\end{lemm}
\begin{pf}
We prove this lemma by reducing the $(2,d)$-defective coloring
problem to the $(\frac{4}{2},d)$-coloring problem. Let $G$ be an
instance of the $(2,d)$-defective coloring problem.

If $G$ has a $(\frac{4}{2},d)$-coloring $f$ using colors in
$\{0,1,2,3\}$, then we define a $(2,d)$-defective coloring $g$ of
$G$ as follows. For each vertex $v\in V(G)$, if $f(v)\in \{0,1\}$
(resp. $f(v)\in \{2,3\}$), then we let $g(v)=0$ (resp. $g(v)=1$). By
the definition of $(\frac{4}{2},d)$-coloring, for any vertex $v\in
V(G)$, if $f(v)\in \{0,1\}$ (resp. $\{2,3\}$), then there are at
most $d$ vertices adjacent to $v$ which are colored by colors in
$\{0,1\}$ (resp. $\{2,3\}$). It follows that under the coloring $g$,
there are at most $d$ vertices adjacent to $v$ which are colored by
the same color as $v$. Thus, $g$ is a $(2,d)$-defective coloring of
$G$.

On the other hand, we suppose $G$ has a $(2,d)$-defective coloring
$g$ using colors in $\{0,1\}$. Then we define a
$(\frac{4}{2},d)$-coloring $h$ of $G$ as follows. For each vertex
$v\in V(G)$, if $g(v)=0$ (resp. $g(v)=1$), then we let $h(v)=0$
(resp. $h(v)=2$). Since for any vertex $v$, there are at most $d$
vertices which have the same color as $v$, it is clear that $h$ is a
$(\frac{4}{2},d)$-coloring of $G$.

We have shown that $G$ has a $(2,d)$-defective coloring if and only
if $G$ has a $(\frac{4}{2},d)$-coloring. Thus the NP-completeness of
the $(\frac{4}{2},d)$-coloring problem follows from that of the
$(2,d)$-defective coloring problem. Therefore, the lemma holds by
Lemma \ref{kt-complexity}.
\end{pf}

At the end of this section, we prove the NP-completeness of Problem
\ref{dis-def-cir-pro} for $k\geq5$. We need the following lemma.

\begin{lemm} \label{k-d-col}
Let $k$ and $d$ be two fixed integers with $k\geq5$ and $d\geq1$.
Then for the complete graph $G=K_{\lfloor\frac{k}{2}\rfloor(d+1)}$,
each vertex in $V(G)$ must be relaxed $d$ times in any
$(\frac{k}{2},d)$-coloring of $G$.
\end{lemm}
\begin{pf}
Suppose $f$ is a $(\frac{k}{2},d)$-coloring of $G$. Suppose to the
contrary that there is a vertex $u\in V(G)$ which has at most $d-1$
relaxations. By the symmetry of colors, we may assume $f(u)=0$.
Because $f$ is a $(\frac{k}{2},d)$-coloring, there are at most $d$
vertices which are assigned colors in $\{0,1,k-1\}$. Therefore, the
number of vertices which are assigned colors in
$A=\{2,3,\dots,k-2\}$ is at most
$\lfloor\frac{k}{2}\rfloor(d+1)-d=(\lfloor\frac{k}{2}\rfloor-1)(d+1)+1$.
Note that $|A|=k-3$. Under the definition of
$(\frac{k}{2},d)$-coloring, we know that two consecutive colors in
$A$ can be assigned to at most $d+1$ vertices. Consequently, if $k$
is even, then the $k-3$ colors in $A$ can be assigned to at most
$\frac{k-4}{2}(d+1)+d+1=\frac{k-2}{2}(d+1)=(\lfloor\frac{k}{2}\rfloor-1)(d+1)$
vertices; if $k$ is odd, then the $k-3$ colors in $A$ can be
assigned to at most
$\frac{k-3}{2}(d+1)=(\lfloor\frac{k}{2}\rfloor-1)(d+1)$ vertices.
This is a contradiction since
$d+(\lfloor\frac{k}{2}\rfloor-1)(d+1)<\lfloor\frac{k}{2}\rfloor(d+1)$.
Thus the lemma holds.
\end{pf}

\begin{lemm} \label{k-d-col-com}
Let $k$ and $d$ be two fixed integers with $k\geq5$ and $d\geq1$.
Then the $(\frac{k}{2},d)$-coloring problem is NP-complete.
\end{lemm}
\begin{pf}
We prove this lemma by reducing the $\frac{k}{2}$-coloring problem
to the $(\frac{k}{2},d)$-coloring problem. Let $G$ be an instance of
the $\frac{k}{2}$-coloring problem. We construct from $G$ an
instance $G^*$ of the $(\frac{k}{2},d)$-coloring problem as follows:
for each vertex $v$ of $G$, take a copy of
$K_{\lfloor\frac{k}{2}\rfloor(d+1)}$ and identify one of its
vertices with $v$. We denote this copy of
$K_{\lfloor\frac{k}{2}\rfloor(d+1)}$ by $K_v$.

Suppose $G$ has a $\frac{k}{2}$-coloring $f$ using colors in
$\{0,1,\dots,k-1\}$. Then we define a $(\frac{k}{2},d)$-coloring of
$G^*$ as follows. For each vertex $v\in V(G)$, let $g(v)=f(v)$. If
$k$ is even (resp. odd), then let the
$\lfloor\frac{k}{2}\rfloor(d+1)$ vertices in $K_v$ be assigned
colors in $\{f(v),f(v)+2,\dots,f(v)+k-2\}$ (resp.
$\{f(v),f(v)+2,\dots,f(v)+k-3\}$) (where all $``+"$s are taken
modulo $k$), such that each color is assigned to $d+1$ vertices. It
is straightforward to check that the coloring $g$ is a
$(\frac{k}{2},d)$-coloring of graph $G^*$.

On the other hand, suppose $G^*$ has a $(\frac{k}{2},d)$-coloring
$g$. Then let $h(v)=g(v)$ for any vertex $v\in V(G)$. By lemma
\ref{k-d-col}, we know that vertex $v$ has $d$ relaxations in
$V(K_v\setminus\{v\})$. Thus the mapping $h$ is a
$\frac{k}{2}$-coloring of $G$.

We have shown that $G$ has a $\frac{k}{2}$-coloring if and only in
$G^*$ has a $(\frac{k}{2},d)$-coloring. So, the NP-completeness of
the $(\frac{k}{2},d)$-coloring is established by Lemma
\ref{comp-k-cir}.
\end{pf}

Combining Lemma \ref{def-k-3}, \ref{k-four-d} and \ref{k-d-col-com},
we have the following result:

\begin{theo} \label{d-def-comp}
Let $k$ and $d$ be two fixed positive integers. The
$(\frac{k}{2},d)$-coloring problem is NP-complete expect the case
$k\leq3$, which is polynomially solvable.
\end{theo}

\section{Relaxed 2-distant circular coloring of outerplanar graphs}

In this section, we study the relaxed 2-distant circular coloring of
outerplanar graphs. An \textit{outerplanar graph} is a graph with a
planar drawing for which all vertices belong to the outer face of
the drawing. It is shown that a graph $G$ is an outerplanar graph if
and only if $G$ has no subdivision of complete graph $K_4$ and
complete bipartite graph $K_{2,3}$. For $d$-defective 2-distant
circular coloring, Klostermeyer studied the
$(\frac{5}{2},1)$-colorability of outerplanar graphs in
[\citeonline{K2}].

\subsection{A graph that is $(\frac{4}{2},1)$-colorable but not $(\frac{4}{2},1)^*$-colorable}

To compare the two relaxations of circular coloring, we give an
example to illustrate the difference between the
$(\frac{k}{2},d)$-coloring and the $(\frac{k}{2},t)^*$-coloring in
this subsection. An outerplanar graph $G_5$ is defined as follows.
The vertex set is $V(G_5)=\{x,y_1,y_2,\dots,y_6,u_1,$
$u_2,\cdots,u_5,v_1,v_2,\dots,v_5\}$. And the edge set is the union
of $\{xy_i|i=1,2,\dots,6\}$ and $\{y_iu_i,u_iv_i,$
$v_iy_{i+1}|i=1,2,\dots,5\}$. Note that the boundary of each
interior face of $G_5$ is a 5-cycle. The vertices along the exterior
face of $G_5$ form a cycle $xy_1u_1v_1y_2\cdots y_5u_5v_5y_6x$.

\begin{figure}[h]
\begin{minipage}{7cm}
\centering \resizebox{4.5cm}{4.5cm}{\includegraphics{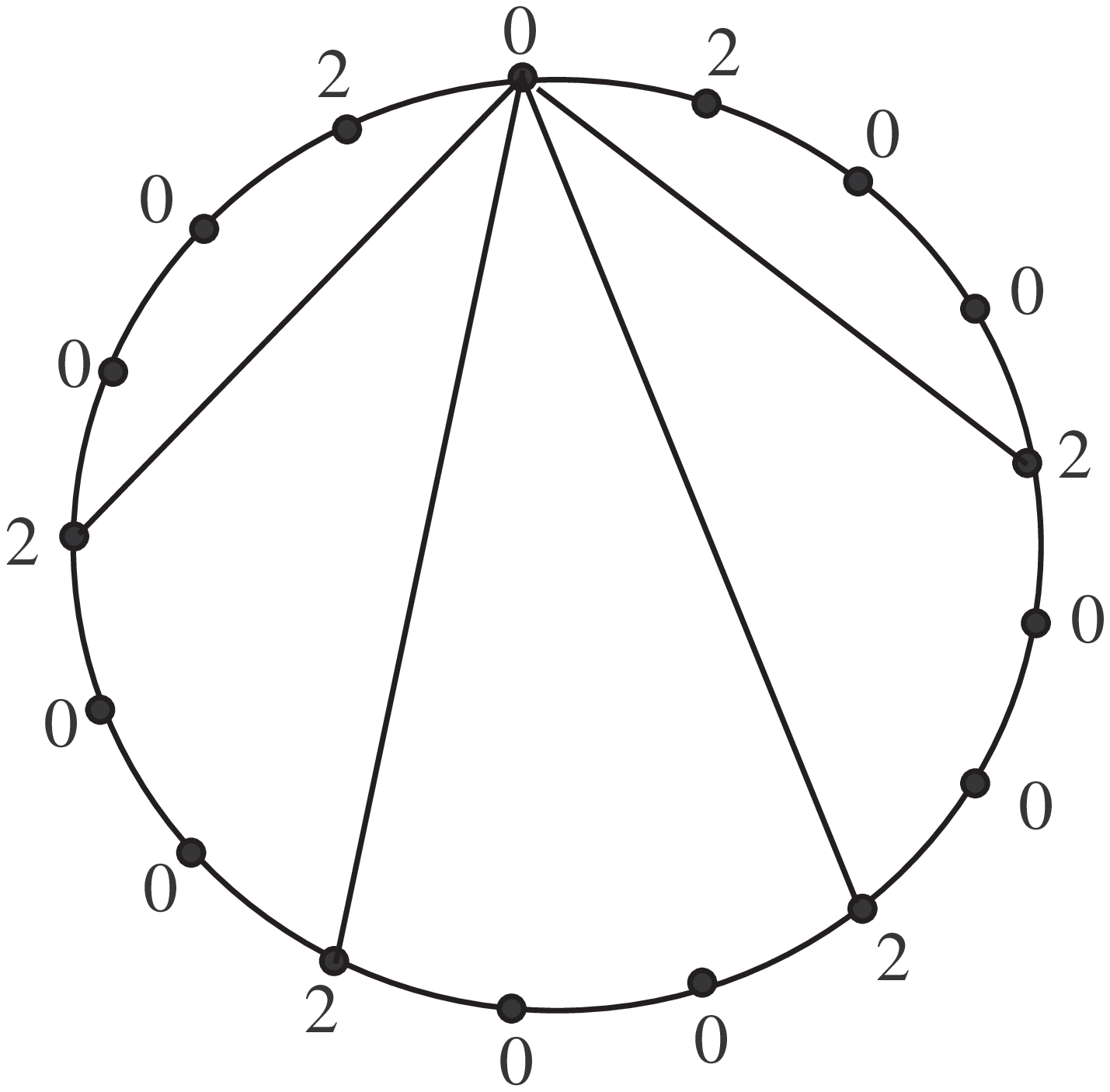}}
\caption{A $(\frac{4}{2},1)$-coloring of $G_5$} \label{dcG5}
\end{minipage}
\hskip 0.6cm
\begin{minipage}{7cm}
\centering\resizebox{4.5cm}{4.5cm}{\includegraphics{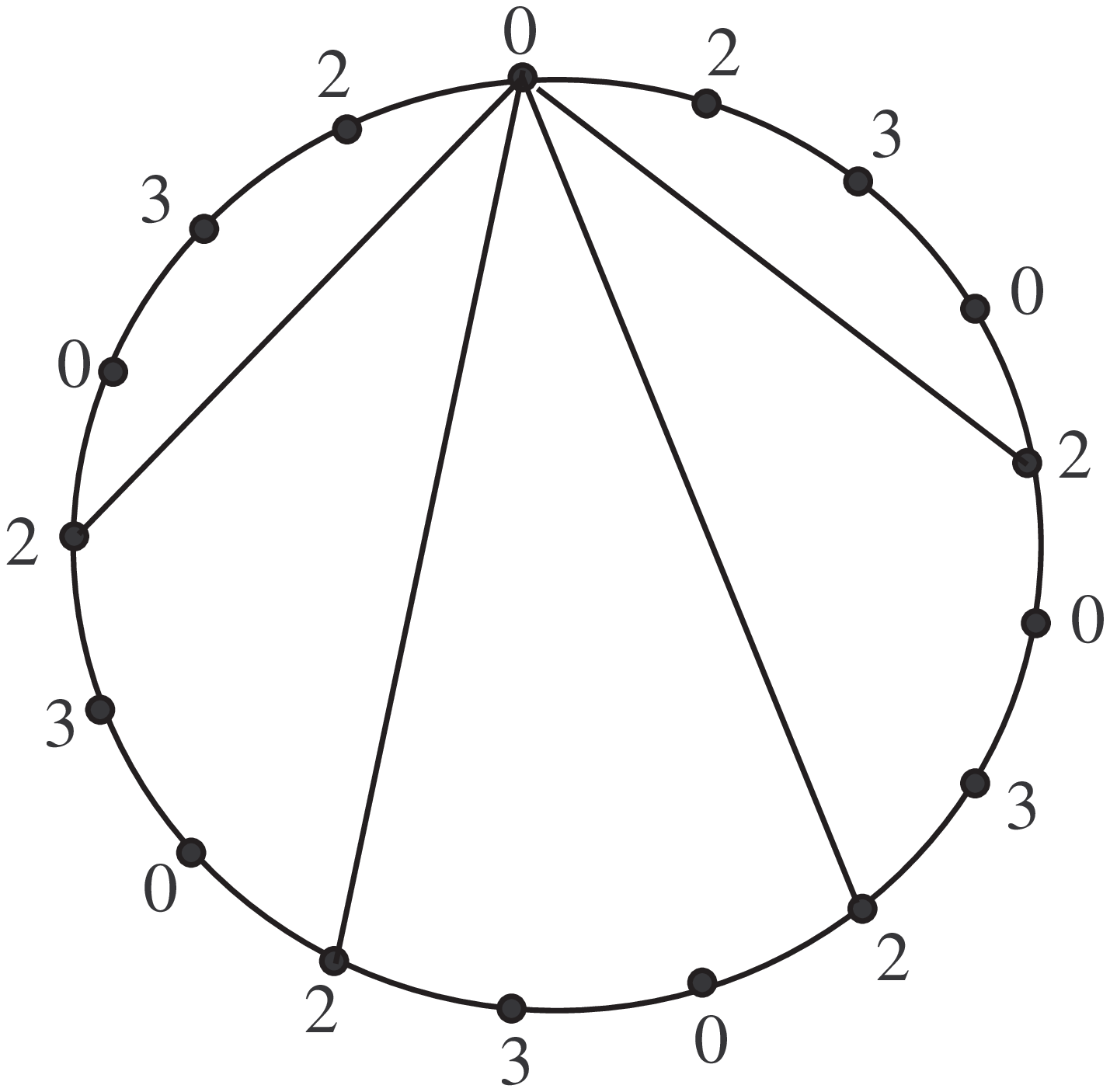}}
\caption{A $(\frac{4}{2},2)^*$-coloring of $G_5$} \label{rcG5}
\end{minipage}
\end{figure}

We First consider the $(\frac{4}{2},1)$-coloring of graph $G_5$. We
give a mapping $f$ from $V(G_5)$ to $\{0,1,2,3\}$ such that
$f(x)=f(u_i)=f(v_i)=0$ ($i=1,2,\dots,5$) and $f(y_j)=2$
($j=1,2,\dots,6$). See Figure \ref{dcG5}. It is clear that the
coloring $f$ satisfies that each vertex $v$ is adjacent to at most
one vertex $u$ with $|f(u)-f(v)|_4\leq1$. Thus, $G_5$ is
$(\frac{4}{2},1)$-colorable. However, for $t$-relaxed 2-distant
circular coloring, we shall prove that the graph $G_5$ is not
$(\frac{4}{2},1)^*$-colorable.

\begin{theo}\label{out-1-relax-no}
The graph $G_5$ is $(\frac{4}{2},2)^*$-colorable, but is not
$(\frac{4}{2},1)^*$-colorable.
\end{theo}

\begin{pf}
The $(\frac{4}{2},2)^*$-colorability of $G_5$ is demonstrated by a
$(\frac{4}{2},2)^*$-coloring of $G_5$ indicated in Figure
$\ref{rcG5}$. We next prove that $G_5$ is not
$(\frac{4}{2},1)^*$-colorable.

Suppose to the contrary that $G_5$ is $(\frac{4}{2},1)^*$-colorable
and let $f: V(G_5)\rightarrow \{0,1,2,3\}$ be a
$(\frac{4}{2},1)^*$-coloring of $G_5$. Without loss of generality,
we may assume $f(x)=0$. We claim that the vertices $y_1,y_2$ and
$y_3$ can not be all colored by $2$. Suppose
$f(y_1)=f(y_2)=f(y_3)=2$. Since $f$ is a $1$-relaxed 2-distant
circular coloring of $G_5$, it is not difficult to see that
$\{f(u_1),f(v_1)\}$ and $\{f(u_2),f(v_2)\}$ must be $\{1,3\}$. But
then $y_2$ is relaxed two times. So, $y_1,y_2$ and $y_3$ can not be
all colored by $2$. It follows that the vertex $x$ must have one
relaxation in $\{y_1,y_2,y_3\}$. Symmetrically, $x$ must also have
one relaxation in $\{y_4,y_5,y_6\}$. This is a contradiction since
$f$ is a $(\frac{4}{2},1)^*$-coloring. Thus the theorem holds.
\end{pf}

\subsection{$(\frac{5}{2},4)^*$-coloring for all outerplanar graphs}

For an outerplanar graph $G$, since the chromatic number of $G$ is
at most 3, we have $c\chi^t_2(G)\leq 2\chi(G)=6$ for any nonnegative
integer $t$, indicating that all outerplanar graphs are
$(\frac{6}{2},0)^*$-colorable. It was proved in [\citeonline{C-C-W}]
that every outerplanar graph is $(2,2)$-defective colorable. With
this result in mind, it is natural to ask whether there exists an
positive integer $t$ such that every outerplanar graph is
$(\frac{5}{2},t)^*$-colorable (resp. $(\frac{5}{2},t)$-colorable).
To answer this question, we first give an outerplanar graph which is
not $(\frac{5}{2},1)^*$-colorable and also is not
$(\frac{5}{2},1)$-colorable (see [\citeonline{K2}]).

\begin{figure}[h!]
\centering \resizebox{7cm}{2.5cm}{\includegraphics{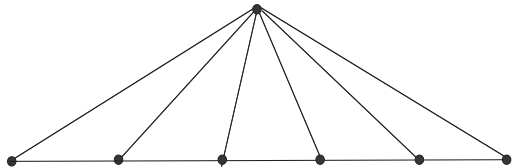}}
\caption{An outerplanar graph which is not
$(\frac{5}{2},1)^*$-colorable} \label{not51}
\end{figure}

\begin{theo}\label{out-1-relax-no}
There exists an outerplanar graph that is not
$(\frac{5}{2},1)^*$-colorable.
\end{theo}
\begin{pf}
It is easy to verify that the graph in Figure \ref{not51} is not
$(\frac{5}{2},1)^*$-colorable. \end{pf}

Next we apply the ordered breadth first search algorithm (which was
introduced by Calamoneri and Petreschi in [\citeonline{C-P}]) for an
outerplanar graph $G$ to construct a spanning tree $T$ of $G$. Then
we produce a particular $(\frac{5}{2},0)^*$-coloring $f$ of $T$.
With the properties of $T$, we finally prove that $f$ is a
$(\frac{5}{2},4)^*$-coloring of the outerplanar graph $G$.

\begin{figure}[h]
\begin{minipage}{7cm}
\centering \resizebox{5cm}{5cm}{\includegraphics{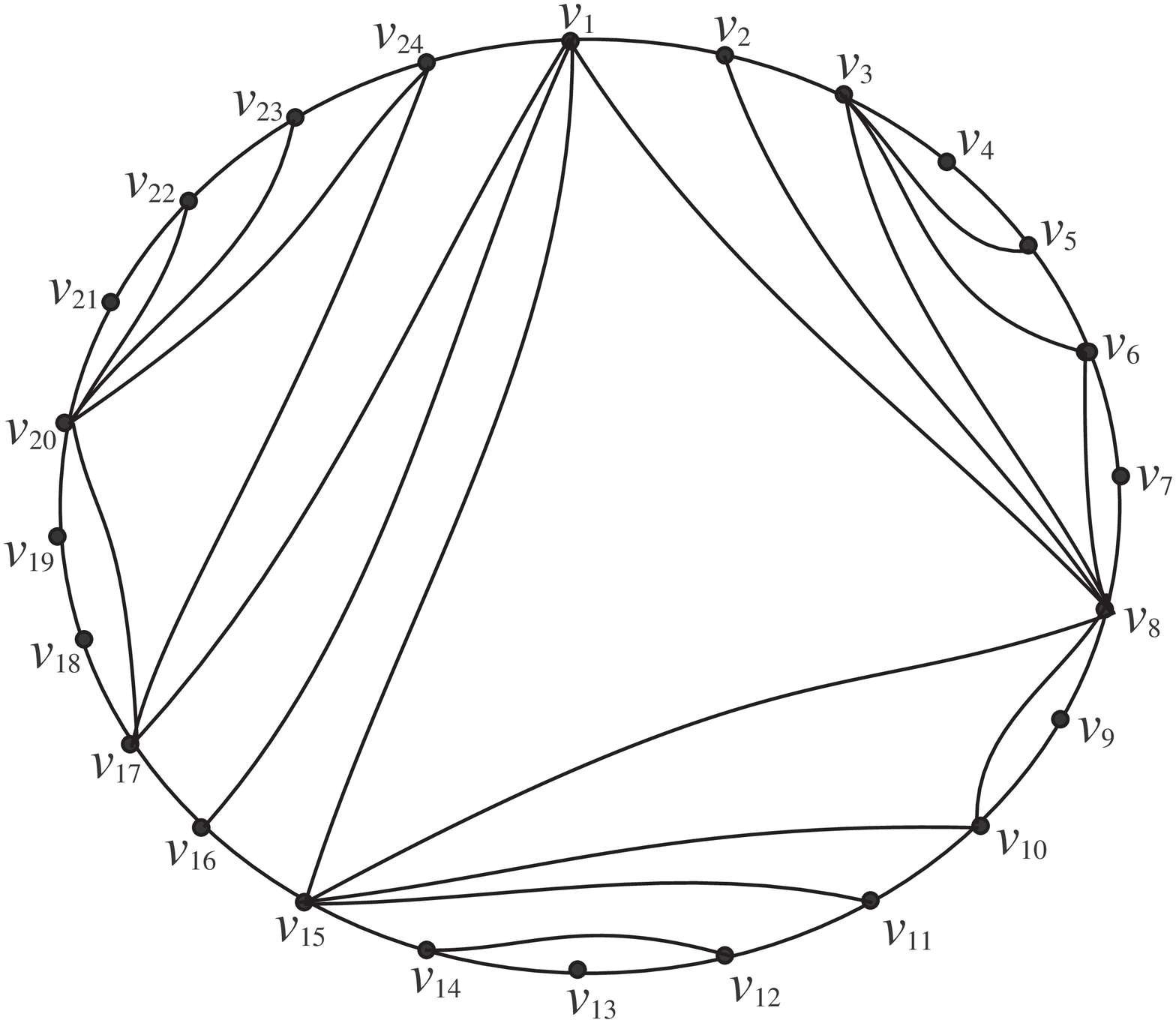}}
\caption{An outerplanar graph $G_1$} \label{G1}
\end{minipage}
\hskip 0.6cm
\begin{minipage}{7.5cm}
\centering\resizebox{6.8cm}{4.8cm}{\includegraphics{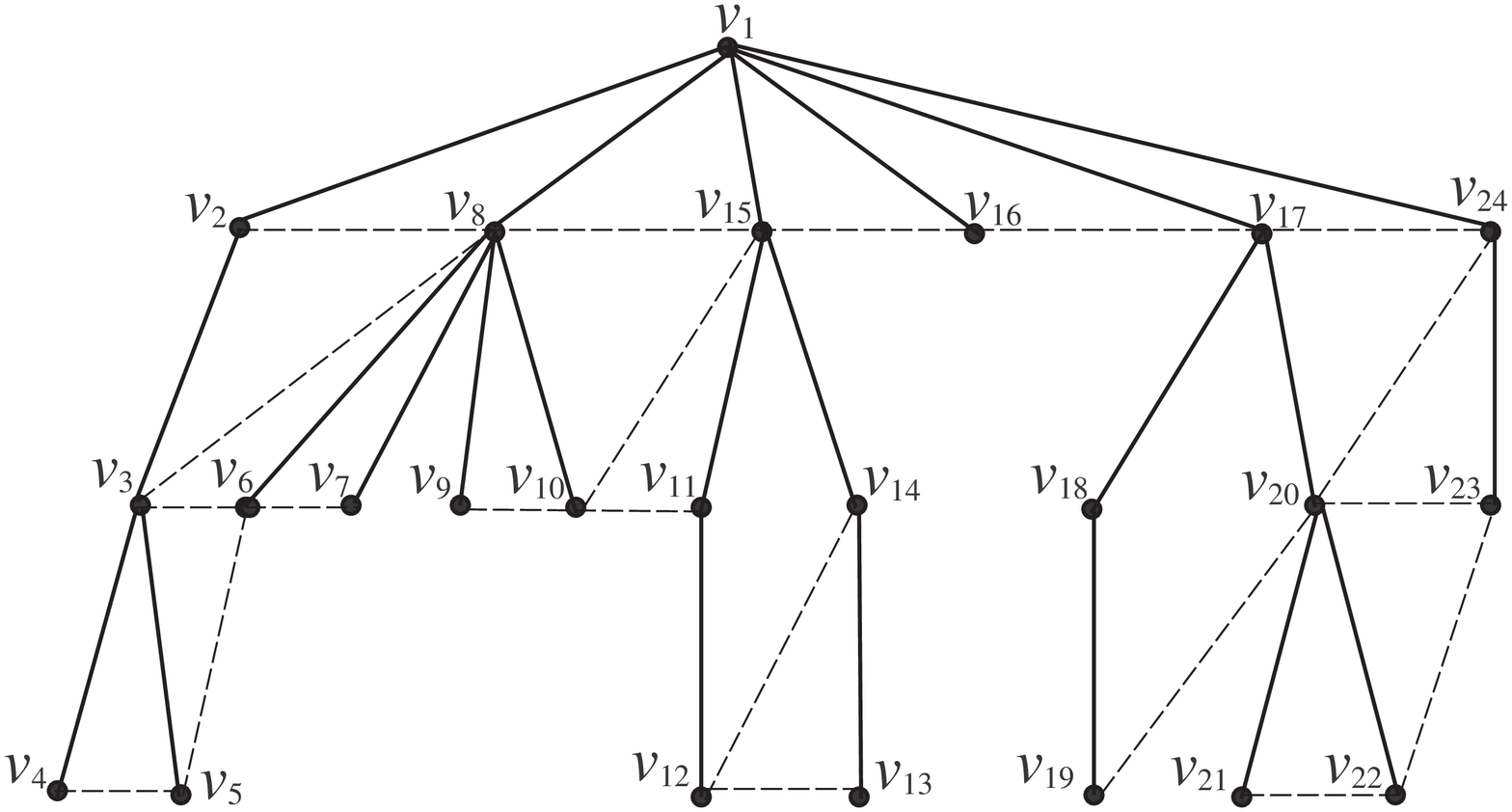}}
\caption{An OBFT partition of $G_1$} \label{G1_OBFT}
\end{minipage}
\end{figure}

Let $G$ be a connected outerplanar graph. Choose a vertex $r\in
V(G)$. At the beginning, we order all vertices clockwise started
from $r$ and along the outerface of $G$. Now, we perform a breadth
first search starting from $r$ in such a way that vertices coming
first in the ordering are visited first. The authors of
[\citeonline{C-P}] called such progress \textit{ordered breadth
first search} (OBFS). The spanning tree obtained by OBFS is called
\textit{ordered breadth first tree} (OBFT). Using OBFS, an connected
outerplanar graph $G$ can be edge-partitioned into a spanning tree
$T$ and a subgraph $H$, that is $E(G)=E(T)\cup E(H)$ and $E(T)\cap
E(H)=\emptyset$. This edge-partition is called an \textit{OBFT
partition}. To give an example of OBFT partition, we consider the
outerplanar graph $G_1$ depicted in Figure \ref{G1}. In Figure
\ref{G1_OBFT}, vertex $v_1$ is the root of the tree $T$ produced by
OBFT partition and solid (resp. broken) lines denote tree-edges
(resp. non-tree-edges). Calamoneri and Petreschi studied the
$L(h,1)$-labeling of outerplanar graphs on the basis of OBFS
algorithm in [\citeonline{C-P}]. Later, the authors of
[\citeonline{W-D-W}] used OBFT partition to design a polynomial-time
algorithm producing nearly-optimal solutions for an edge coloring
problem of outerplanar graphs.

A tree with a particular vertex $r$ designated as its root is called
\textit{a rooted tree}, denoted by $T_r$. According to the distance
to the root $r$, the vertices of a rooted tree can be arranged in
layers. All vertices at distance $i$ from the root form the
$i^{th}$-layer. Hence, the $0^{th}$-layer consists of the root only.
For a vertex $u$ in the $i^{th}$-layer ($i\geq1$), the neighbor of
$u$ in the $(i-1)^{th}$-layer is called its \textit{father} and all
neighbors of $u$ in the $(i+1)^{th}$-layer are called its
\textit{sons}. Let $SON(u)$ denote the set of all sons of $u$. If
two vertices have the same father, we say they are
\textit{brothers}. Given two vertices $u$ and $v$ in $V(T_r)$, let
$P(u,v)$ denote the unique path from $u$ to $v$ in $T_r$. The last
common vertex of $P(r,u)$ and $P(r,v)$ is called \textit{the latest
common ancestor} of $u$ and $v$, denoted by $LCA(u,v)$. Let $u$ and
$v$ be two vertices of $G$. If $uv\in E(H)$ or $u$ and $v$ are two
consecutive vertices in the same layer, then we define $Int(u,v)$ as
the set of vertices falling within the region bounded by two paths
$P(LCA(u,v),u)$ and $P(LCA(u,v),v)$ together with the edge $uv$. For
example, in Figure \ref{G1_OBFT}, we have
$Int(v_{10},v_{11})=\emptyset$, $Int(v_{13},v_{19})=\{v_{16}\}$ and
$Int(v_{19},v_{20})=\emptyset$. For a connected outerplanar graph,
choosing any vertex $r$ and using OBFT partition, we can get a
rooted tree $T_r$ and a subgraph $H$. Let $L(T_r)$ denote the number
of layers of $T_r$. For $i=0,1,\dots, L(T_r)-1$, denote by $l_i$ the
number of vertices in the $i^{th}$-layer. Vertices in the
$i^{th}$-layer are denoted from left to right by
$v_1^i,v_2^i,\dots,v_{l_i}^i$. Some key properties about $T_r$ and
$H$ are either given in [\citeonline{C-P}] and [\citeonline{W-D-W}]
or easily obtained according to outerplanarity of graph.

\begin{lemm}\label{OBFT-pro}
Every OBFT partition $T_r\cup H$ for a connected outerplanar graph
$G$ has the following properties:

\noindent (1) The maximum degree of $H$ is at most $4$.

\noindent (2) If $v_p^i$ is adjacent to $v_q^i$, then $v_p^iv_q^i$
is a non-tree-edge, and $q=p+1$ or $p=q+1$.

\noindent (3) If $v_p^iv_q^j\in E(H)$ and $j<i$, then $j=i-1$.

\noindent (4) If $v_p^{i}v_q^{i-1}\in E(H)$ and $v_p^i$ is a son of
$v_h^{i-1}$, then $q=h+1$ and $v_p^i$ is the rightmost son of
$v_h^{i-1}$. Meanwhile, $Int(v_p^{i},v_q^{i-1})=\emptyset$.

\noindent (5) Let $v_p^{i+1}$ be a son of $v_l^i$, and $v_q^{i+1}$
be a son of $v_h^i$ such that $l\leq h$ and $p<q$. If
$v_p^{i+1}v_q^{i+1}\in E(H)$, then $q=p+1$, and $h=l$ or $h=l+1$.
Meanwhile, $Int(v_p^{i+1},v_q^{i+1})=\emptyset$.
\end{lemm}

\begin{figure}[h!]
\centering
\resizebox{8.5cm}{4.5cm}{\includegraphics{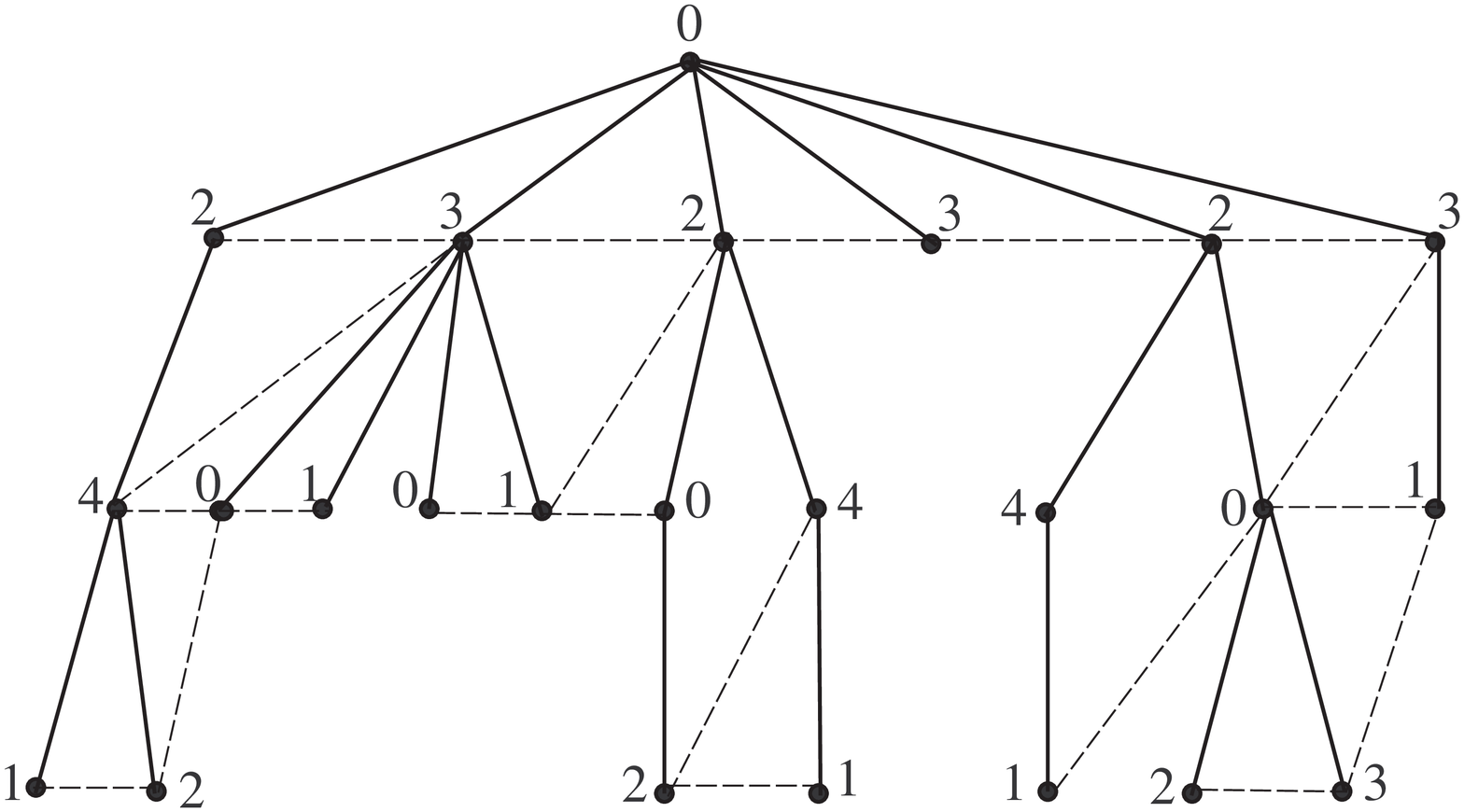}}
\caption{A $(\frac{5}{2},0)^*$-coloring for subgraph $T_{v_1}$ of
$G_1$ constructed by Algorithm 1} \label{G1coloring}
\end{figure}

For a connected outerplanar graph $G$, we choose any vertex $r$ as
root and use OBFT partition to get a rooted tree $T_r$. The
following algorithm will produce a $(\frac{5}{2},0)^*$-coloring of
$T_r$ and Figure \ref{G1coloring} gives a
$(\frac{5}{2},0)^*$-coloring constructed by the algorithm for
subgraph $T_{v_1}$ of $G_1$.

\begin{algorithm}[h!]
\caption{Construct a $(\frac{5}{2},0)^*$-coloring $f$ for $T_r$}
\begin{algorithmic}[1]
\STATE Assign color $0$ to the root $r$, that is $f(r)=0$; \\
\STATE Assign colors $2$ and $3$ to the sons of $r$ alternately from
left to right;\\
 \FOR{$i:=1$ to $L(T_r)-2$}
       \STATE Let $q=\min\limits_{1\leq h\leq
       i}\{SON(v_h^i)\neq\emptyset\}$;\\
       Assign colors $f(v_q^i)+2 \pmod 5$ and $f(v_q^i)+3 \pmod 5$ to the sons of $v_q^i$
alternately from left to right;\\
        \FOR{$p:=q+1$ to $l_i$}
   \WHILE{$SON(v_p^i)\neq\emptyset$}
      \IF{$SON(v_{p-1}^i)\neq\emptyset$ and $Int(v_{p-1}^i,v_p^i)=\emptyset$}
        \STATE Assign colors $f(v_p^i)+2 \pmod 5$ and $f(v_p^i)+3 \pmod 5$ to the sons of $v_p^i$ alternately
  from left to right, such that the circular distance between the
color assigned to the leftmost son of $v_p^i$ and the color assigned
to the rightmost son of $v_{p-1}^i$ is one.
       \ELSE
          \STATE Assign colors $f(v_p^i)+2 \pmod 5$ and $f(v_p^i)+3 \pmod 5$ to the sons of $v_p^i$ alternately from
left to right.
    \ENDIF \ENDWHILE
 \ENDFOR \ENDFOR
\end{algorithmic}\label{52-c-Tr}
\end{algorithm}

To prove the correctness of the algorithm, it suffices to verify
that in line 8 of Algorithm 1, sons of $v_p^i$ can be assigned the
colors satisfying the conditions specified in the algorithm.

\begin{lemm}\label{correctness}
Let $G$ be an outerplanar graph and $T_r$ a rooted tree obtained by
an OBFT partition of $G$. For $i\in \{1,2,\dots, L(r)-2\}$ and $p\in
\{1,2,\dots,l_i-1\}$, let $v_h^{i+1}$ be the rightmost son of
$v_{p}^i$ and $v_{h+1}^{i+1}$ be the leftmost son of $v_{p+1}^i$.
Suppose $f$ is the coloring of $T_r$ produced by Algorithm 1. If
$Int(v_{p}^i,v_{p+1}^i)=\emptyset$, then
$|f(v_h^{i+1})-f(v_{h+1}^{i+1})|_5=1$.
\end{lemm}

\begin{pf}
We prove this lemma by induction on the value of $i$. Assume $i=1$.
By Algorithm 1, we have $\{f(v_{p}^1),f(v_{p+1}^1)\}=\{2,3\}$ for
all $p\in\{1,2,\dots,l_1-1\}$. Without loss of generality, we assume
$f(v_{p}^1)=2$ and $f(v_{p+1}^1)=3$. According to Algorithm 1,
colors $0$ and $4$ appear on the sons of $v_{p}^1$ alternately from
left to right. If $f(v_h^{2})=0$ (resp. $f(v_h^2)=4$), then
Algorithm 1 will assign color 1 (resp. color 0) to $v_{h+1}^2$, and
all sons of $v_{p+1}^1$ from left to right are colored in the
pattern as $(1010\dots)$ (resp. $(0101\dots)$). It is clear that the
lemma holds for $i=1$.

Now let $i$ be an integer in $\{2,3,\dots,L(r)-2\}$. Let
$LCA(v_{p}^i,v_{p+1}^i)=v_{h_1}^{i-j}$ for some $1\leq j\leq i$ and
$1\leq h_1\leq l_{i-j}$. Let
$P(v_{h_1}^{i-j},v_{p}^i)=v_{h_1}^{i-j}v_{h_2}^{i-j+1}\dots
v_{h_{j}}^{i-1}v_{p}^i$. Because $Int(v_{p}^i,v_{p+1}^i)=\emptyset$,
$P(v_{h_1}^{i-j},v_{p+1}^i)$ must be
$v_{h_1}^{i-j}v_{h_2+1}^{i-j+1}\dots v_{h_{j}+1}^{i-1}v_{p+1}^i$. If
$j=1$, then $v_{p}^i$ and $v_{p+1}^i$ are brothers and so according
to Algorithm 1, we have $|f(v_{p}^i)-f(v_{p+1}^i)|_5=1$. If $j\geq
2$, since $Int(v_{h_{j}}^{i-1},v_{h_{j}+1}^{i-1})=\emptyset$, by the
inductive hypothesis, we also have $|f(v_{p}^i)-f(v_{p+1}^i)|_5=1$.
Without loss generality, we may assume $f(v_{p}^i)=a$ and
$f(v_{p+1}^i)=a+1\pmod5$ for some $a\in\{0,1,2,3,4\}$. Then
Algorithm 1 assigns $a+2\pmod5$  or $a+3\pmod5$ (resp. $a+3\pmod5$
or $a+4\pmod5$) to $v_h^{i+1}$ (resp. $v_{h+1}^{i+1}$). According to
Algorithm 1, if $f(v_h^{i+1})=a+2\pmod5$ (resp.
$f(v_h^{i+1})=a+3\pmod5$), then $f(v_{h+1}^{i+1})=a+3\pmod5$ (resp.
$f(v_{h+1}^{i+1})=a+4\pmod5$. In both two cases, we have
$|f(v_h^{i+1})-f(v_{h+1}^{i+1})|_5=1$. The proof is completed.
\end{pf}

Note that the coloring $f$ produced by Algorithm 1 satisfies
$|f(u)-f(v)|_5=2$ for any edge $uv\in E(T_r)$. The following result
is obtained.

\begin{lemm}\label{Tr-c-coloring}
Suppose $G$ is a connected outerplanar graph, Algorithm 1 produces a
$(\frac{5}{2},0)^*$-coloring for any ordered breadth first tree of
$G$.
\end{lemm}

We next show that the $(\frac{5}{2},0)^*$-coloring $f$ of $T_r$
produced by Algorithm 1 is actually a $(\frac{5}{2},4)^*$-coloring
of $G$. Thus we obtain the following theorem.

\begin{theo} \label{524star} All Outerplanar graphs are
$(\frac{5}{2},4)^*$-colorable.
\end{theo}
\begin{pf}
Let $G$ be any outerplanar graph $G$. Without loss of generality,
assume $G$ is connected. Choose any vertex $r\in V(G)$ as root, we
use OBFT partition to get two edge disjoint subgraphs $T_r$ and $H$
of $G$. Note that $V(T_r)=V(G)$. Let $f$ be the
$(\frac{5}{2},0)^*$-coloring $f$ of $T_r$ produced by Algorithm 1.
Then $|f(u)-f(v)|_5=2$ for any edge $uv\in E(T_r)$. Since the
maximum degree of $H$ is at most $4$ by Lemma \ref{OBFT-pro} (1), we
only need to prove that $f(x)\neq f(y)$ for each edge $xy\in E(H)$.

\textbf{Case 1}: $x$ and $y$ are in the same layer. Suppose
$x=v_p^i$. Since $xy\in E(H)$, we may assume $y=v_{p+1}^i$ by Lemma
\ref{OBFT-pro} (2). Let $v_q^{i-1}$ be the father of $v_p^i$ and
$v_h^{i-1}$ the father of $v_{p+1}^i$. From Lemma \ref{OBFT-pro}
(5), we have $h=q$ or $h=q+1$. If $h=q$, then $v_p^i$ and
$v_{p+1}^i$ are brothers. It is obvious that $f(v_p^i)\neq
f(v_{p+1}^i)$ by Algorithm 1. If $h=q+1$, then $v_p^i$ is the
rightmost son of $v_q^{i-1}$ and $v_{p+1}^i$ is the leftmost son of
$v_h^{i-1}$. Because $v_p^i v_{p+1}^i\in E(G)$, we have
$Int(v_p^i,v_{p+1}^i)=\emptyset$ by Lemma \ref{OBFT-pro} (5),
implying $Int(v_q^{i-1},v_h^{i-1})=\emptyset$. Hence, $f(v_p^i)\neq
f(v_{p+1}^i)$ according to Algorithm 1.

\textbf{Case 2}: $x$ and $y$ are in different layers. By Lemma
\ref{OBFT-pro} (4), we may assume that $x$ is the rightmost son of
$v_q^{i-1}$ and $y=v_{q+1}^{i-1}$ (where $i\geq2$) with
$Int(x,v_{q+1}^{i-1})=Int(v_q^{i-1},v_{q+1}^{i-1})=\emptyset$.
According to Algorithm 1, we must have $|f(v_q^{i-1})-f(x)|_5=2$ and
$|f(v_q^{i-1})-f(y)|_5=1$, implying $f(x)\neq f(y)$.
\end{pf}

For $d$-defective 2-distant circular coloring, we can conclude the
following result with the aid of ordered breadth first search
algorithm.

\begin{theo} \label{422} All Outerplanar graphs are $(\frac{4}{2},2)$-colorable.
\end{theo}

\begin{pf}
Let $G$ be any outerplanar graph. Without loss of generality, assume
$G$ is connected. Choose any vertex $r$ as root and use OBFS
partition to get a rooted tree $T_r$ and a subgraph $H$. We define a
mapping $f$ from $V(T_r)$ to $\{0,1,2,3\}$ as follows: let $f(r)=0$,
$f(v)=2$ for any vertex $v$ in odd layers, and $f(u)=0$ for any
vertex $u$ in even layers. It is clear that $|f(x)-f(y)|_4=2$ if
$xy$ is an edge of $T_r$ or $xy$ is an edge of $H$ with $x$ and $y$
in different layer, and $f(x)=f(y)$ if $xy$ is an edge of $H$ with
$x$ and $y$ in the same layer. By Lemma \ref{OBFT-pro} (2), any
vertex $u$ of $G$ has at most two neighbors with $|f(u)-f(v)|_4<2$.
Thus $f$ is a $(\frac{4}{2},2)$-coloring of $G$.
\end{pf}

Regarding to the conclusions in Theorem \ref{524star} and Theorem
\ref{422}, it is now natural to ask if there exists a fixed positive
integer $t$ such that every outerplanar graph is
$(\frac{4}{2},t)^*$-colorable. We answer this question in next
subsection.

\subsection{$(\frac{4}{2},t)^*$-colorability of outerplanar graphs}

In this subsection, we consider the $(\frac{4}{2},t)^*$-colorability
of outerplanar graphs. For a nonnegative integer $m$, we define a
graph, denoted by $H(m)$, as follows. The vertex set
$V(H(m))=\{x,x_1,x_2,\dots,x_m\}\cup\{y,y_1,y_2,\dots,y_m\}\cup\{z,z_1,z_2,\dots,z_m\}$.
And the edge set is the union of $\{x_ix_{i+1}|i=1,2,\dots,m-1\}$,
$\{y_iy_{i+1}|i=1,2,\dots,m-1\}$, $\{z_iz_{i+1}|i=1,2,\dots,m-1\}$,
$\{x_my,y_mz,z_mx,xy,yz,zx\}$, $\{xx_i|i=1,2,\dots,m\}$,
$\{yy_i|i=1,2,\dots,m\}$ and $\{zz_i|i=1,2,\dots,m\}$. It is obvious
that $H(m)$ is an outerplanar graph and $xx_1\dots x_myy_1\dots
y_mzz_1\dots z_m$ are in the exterior face in counterclockwise
direction. Let $RT_f(x)$ be the number of relaxations of vertex $x$
under the coloring $f$.

\begin{lemm}\label{h2t-lemma}
Let $t$ be a positive integer and $H(2t)$ the graph defined above.
Then for any relaxed 2-distant circular 4-coloring $f$ of $H(2t)$,
there exists at least one vertex in $\{x,y,z\}$ which has at least
$t+2$ relaxations.
\end{lemm}

\begin{pf} Let $f$ be a relaxed 2-distant circular 4-coloring of
$H(2t)$. Suppose to the contrary that $RT_f(x)$, $RT_f(y)$ and
$RT_f(z)$ are less than or equal to $t+1$. Without loss of
generality, assume $f(x)=0$. Note that
$N(x)=\{x_1,x_2,\dots,x_{2t},y, z,z_{2t}\}$ and $x_1x_2,
x_2x_3,\dots,x_{2t-1}x_{2t},x_{2t}y,yz$ and $zz_{2t}$ are edges of
$H(2t)$, at most $t+2$ neighbors of $x$ are assigned color $2$ and
at least $t+1$ neighbors of $x$ are assigned colors in $\{1,3\}$.
Thus $RT_f(x)\geq t+1$ and so $RT_f(x)=t+1$. It follows that
$f(y)=1$ and $f(z)\in\{1,3\}$. Similarly, we have $RT_f(y)=t+1$. By
the symmetry of colors, $RT_f(y)=t+1$ together with $f(y)=2$ implies
$f(z)=0$. This is a contradiction. The lemma follows.
\end{pf}


Note that the graph $H(0)$ is a 3-cycle with vertices $x,y$ and $z$.
It is easy to see that $H(0)$ is not $(\frac{4}{2},1)^*$-colorable.
For any integer $t$ with $t\geq2$, we have the following result:

\begin{theo}\label{h2t-2}
Let $t\geq2$ be a positive integer and $H(2t-2)$ the graph defined
above. Then $H(2t-2)$ is $(\frac{4}{2},t+1)^*$-colorable, but is not
$(\frac{4}{2},t)^*$-colorable.
\end{theo}

\begin{pf}
To prove that $H(2t-2)$ is $(\frac{4}{2},t+1)^*$-colorable, it
suffices to give a $(\frac{4}{2},t+1)^*$-coloring of $H(2t-2)$. We
define a mapping $f$ from $V(H(2t-2))$ to $\{0,1,2,3\}$ as follows.
Let $f(x)=0, f(y)=1$ and $f(z)=2$. For $i\in \{1,2,\dots,t-1\}$, let

$$\left\{\begin{array}{lll} f(x_{2i-1})=2,\\[1mm]
f(x_{2i})=f(y_{2i-1})=3,\\[1mm]
f(y_{2i})=f(z_{2i-1})=0,\\[1mm]
f(z_{2i})=1.
\end{array}\right.$$

\noindent It is clear that $RT_f(x)=RT_f(y)=t+1$ and $RT_f(z)=t$.
For each vertex $u\in V(H(2t-2))\setminus\{x,y,z\}$, as $d(u)\leq 3$
and $t\geq2$, $RT_f(u)\leq t+1$. Therefore, $f$ is a
$(\frac{4}{2},t+1)^*$-coloring of $H(2t-2)$. We illustrate a
$(\frac{4}{2},t+1)^*$-coloring $f$ of $H(6)$ ($t=4$) in Figure
\ref{h6}.

\begin{figure}[h!]
\centering \resizebox{4.5cm}{4.5cm}{\includegraphics{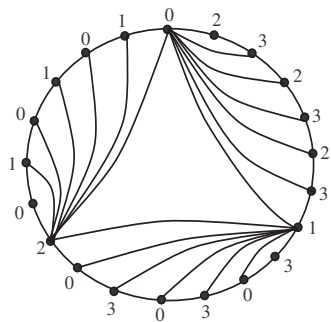}}
\caption{A $(\frac{4}{2},5)^*$-coloring of $H(6)$}\label{h6}
\end{figure}

By Lemma \ref{h2t-lemma}, we know that for any relaxed 2-distant
circular 4-coloring of $H(2t-2)$, there exists at least one vertex
in $\{x,y,z\}$ which has at least $t+1$ relaxations. That is to say,
$H(2t-2)$ can not be $(\frac{4}{2},t)^*$-colorable and the theorem
holds.
\end{pf}

Note that Theorem \ref{h2t-2} shows there does not exist a fixed
positive integer $t$ such that all outerplanar graphs are
$(\frac{4}{2},t)^*$-colorable.

We end the paper by proposing some further research problems on
$(\frac{k}{2},d)$-coloring and the $(\frac{k}{2},t)^*$-coloring of
graphs.

\noindent \textmd{\textbf{Problem 1}}: We have shown that every
outerplanar graph is $(\frac{5}{2},4)^*$-colorable and there is an
outerplanar graph which is not $(\frac{5}{2},1)^*$-colorable. It
will be interesting to determine the smallest integer $t$ such that
every outerplanar graph is $(\frac{5}{2},t)^*$-colorable.

\noindent \textmd{\textbf{Problem 2}}: Is there a fixed integer $t$
such that every planar graph is $(\frac{7}{2},t)^*$-colorable?

\noindent \textmd{\textbf{Problem 3}}: It is meaningful to clarify
the complexity of the $(\frac{k}{2},d)$-coloring problem (resp. the
$(\frac{k}{2},t)^*$-coloring) in planar graphs.

\noindent \textmd{\textbf{Problem 4}}: We expect sufficient
conditions under which every outerplanar graph is
$(\frac{5}{2},1)^*$-colorable and sufficient conditions under which
every outerplanar graph is $(\frac{4}{2},t)^*$-colorable for some
fixed integer $t$.

\noindent \textmd{\textbf{Problem 5}}: It is natural to investigate
the $(\frac{k}{q},t)^*$-coloring problem for $q\geq3$.




\end{document}